%% file: arxiv.tex
\documentclass[10pt]{siamltex}
\usepackage{amssymb,latexsym,graphicx, amsmath,multicol}

\newcommand{\polyf}{\mathcal{P}_d}
\newcommand{\polyfe}{\mathcal{P}_{d+e}}

\DeclareMathOperator{\Span}{span}

\begin{document}

\title{
Several new quadrature formulas for polynomial 
integration in the triangle
}

\author{Mark A. Taylor\thanks{Sandia National Laboratory, 
     Albuquerque, NM, USA. mataylo@sandia.gov} \and Beth A. 
Wingate\thanks{Los Alamos National Laboratory, Los Alamos, NM, USA.
     wingate@lanl.gov} \and Len P. Bos \thanks{Department of Mathematics,
     University of Calgary, Calgary, Alberta Canada
     \mbox{lpbos@math.ucalgary.ca} }}

\maketitle
\bibliographystyle{siam}

\begin{abstract}
We present several new quadrature formulas in the triangle
for exact integration of polynomials.
The points were computed numerically 
with a cardinal function
algorithm which imposes that the 
number of quadrature points $N$ be equal to the
dimension of a lower dimensional polynomial space. 
Quadrature forumulas are presented for up to degree $d=25$, all 
which have positive weights and contain
no points outside the triangle.  Seven of these quadrature
formulas improve on previously known results.
\end{abstract}



\section{Introduction}
We consider a set of $N$ points 
$\{ z_1, z_2, \dots, z_N \}$ 
and assocated weights $\{ w_1, w_2, \dots, w_N \}$ to be 
a quadrature formula of strength $d$ 
if the quadrature approximation for a domain $\Omega$,
\[ 
\int_\Omega g \simeq  \sum_{j=1}^{N}  w_j g(z_j),
\] 
is exact for all polynomials $g$ up to degree $d$.  Among all
quadrature formulas of strength $d$, the optimal ones are those
with the fewest possible points $N$. 
The quadrature problem has been extensively studied and 
has a long history of both theoretical
and numerical development.  For a recent review, see 
\cite{Cools1992, LynessCools1994, CoolsRab1993, Cools1999}.  
An on-line database containing the many of the best known quadrature 
formula is described in \cite{Cools2003}.
Much of these results are also collected and distributed on
CD-ROM in the book \cite{Solin03}.

One successful approach for numerically finding quadrature formulas
dates to \cite{Lyness1975}.  A generalized
version was used recently in \cite{WandzuraXiao2003}.  Newton's method
is used to solve the nonlinear system of algebraic equations for the
quadrature weights and locations of the points.  Symmetry is used to
reduce the complexity of the problem.  If the quadrature points are
invariant under the action of a group $G$, then the number of
equations can be reduced to the dimension of the polynomial subspace 
invariant under $G$.

Recently, a cardinal function algorithm has been developed which can
provide additional reduction in the complexity of the quadrature problem
\cite{TaylorWingateBos2005}.  It is motivated by a similar 
cardinal function Fekete point algorithm
\cite{TaylorWingateVincent2000}.  The key idea is to 
looks for quadrature formula
that have the same number of points as the dimension of a 
lower dimensional polynomial space.  One can then 
construct a cardinal function basis for this lower dimensional
space, make use of a multi-variate generalization of the Newton-Cotes 
quadrature weights and derive a remarkable 
expression analytically relating the variation in the quadrature
weights to the variation of the quadrature points.  The net result
is a significant reduction in the number of equations and 
unknowns, while still retaining analytic expressions for the gradients
necessary to apply steepest decent or Newton iterations.

Symmetry can still be used to further reduce
the complexity of the problem if needed.  However here we have
been able to find optimal quadrature sets of strength 9 through 
25, subject only to the  cardinal function constraint without imposing any 
symmetry constraints on the solutions.

\section{Notation}

Let $\xi = ( \xi_1, \xi_2)$ be an arbitrary point in $\Re^2$.
We will work in the right triangle, $\xi_1 \ge -1$,
$\xi_2 \ge -1$ and $\xi_1+\xi_2 \le 0$.  Let 
$\polyf$ be the finite dimensional vector space of polynomials 
of at most degree $d$, 
\[
\polyf = \Span \{ \xi_1^n \xi_2^m , m+n \le d\}.
\]
We also define
\[
N=\dim \polyf = \frac12(d+1)(d+2).
\] 
The monomials $\xi_1^n \xi_2^m$ are notoriously ill-conditioned, so
it is necessary to describe $\polyf$ with a more reasonable
basis.  For this we use the orthogonal 
Kornwinder-Dubiner polynomials $\{ g_{m,n} \}$ \cite{Appell,Koornwinder,Dubiner},  
\[
g_{m,n}(\xi) = P^{0,0}_m \left( \frac{\xi_1}{1-\xi_2} \right) (1-\xi_2)^m 
P^{2m+1,0}_n (\xi_2) 
\]
where $P^{\alpha,\beta}_n$ are the Jacobi Polynomials with weight
$(\alpha,\beta)$ and degree $n$.  In this basis, $\polyf = \{ g_{m,n},  m+n \le d \}$.
Suitable recurrence
relations for these polynomials are given in 
\cite{SherwinKarniadakis1999}.

\section{Quadrature formula with $N$ points}
We note that given a set of $N$ non-degenerate points in
the triangle $\{ z_j \}$, we can obtain the generalized Newton-Cotes
weights by solving the  $N \times N$ system:
\begin{equation}
\label{E:quad_basis}
\sum_{j=1}^{N}  w_j g_{m,n}(z_j) =  \int g_{m,n} \,d\xi   \qquad \forall g_{m,n} \in \polyf
\end{equation}
By construction, the Newton-Cotes
weights and the points $\{ z_i \}$ give a quadrature formula which
exactly integrates our $N$ basis functions, and thus
\[
\sum_{j=1}^N w_j g(z_j) = \int g \,d\xi  \qquad \forall g \in \polyf.
\]
Because any set of $N$ quadrature points of strength $d$ or greater 
must satisfy Eq.~\ref{E:quad_basis}, the weights for 
all such quadrature formulas must be the Newton-Cotes weights.

To obtain quadrature points of strength greater than $d$, one must optimize
the location of the points $\{ z_j \}$.  Here we present results using
the algorithm from \cite{TaylorWingateBos2005} to perform this
optimization.  The goal is to find points $\{ z_j \}$ which 
exactly integrate all of $\polyfe$ for the largest possible $e$.

\section{Degrees of Freedom bound}
There is one degree of freedom for each coordinate of each 
point, for a total of $2 \dim \polyf$.  Since we are using Newton-Cotes
weights, we automatically integrate exactly all of $\polyf$.  The
number of additional equations that must be satisfied to integrate exactly all
of $\polyfe$ is thus $\dim \polyfe - \dim\polyf$.  
If we require that the degrees-of-freedom
in the location of the quadrature points is at least as large
as the number of equations that must be satisfied, we arrive
at a lower bound on $N$ given by $\dim \polyfe \le   3 \dim \polyf = 3 N$.

\section{Results}
\label{section:results}
Our results for the triangle are summarized in Table~\ref{T:quad}.
Except for quadrature formulas associated with $d=3$ and $d=4$, 
were were able to obtain the optimal solution
(fewest number of points) subject to the cardinal function
constraint on the number of points and the degrees-of-freedom
lower bound:
\begin{align}
\label{E:optimal1}
  &N = \dim \polyf, \\
\label{E:optimal2}
 &\dim \polyfe \le 3 N.
\end{align}
All the quadrature points have positive weights and no points lie outside
the triangle, although neither of these properties is in any way
guaranteed by the cardinal function algorithm.  
The errors presented in the table is the max norm
of the quadrature error over all the ortho-normal basis functions:
\[
\max_{g_{m,n} \in \polyfe} 
\left| \sum_i w_i g_{m,n}(z_i) - \int g_{m,n} \,d\xi \right|
\]
with normalization $\int g_{m,n}^2 \,d\xi = 2$ (the area of
the right triangle). 
Many of the quadrature sets are invariant under the symmetry group of
rotations and reflections of the triangle, $D_3$.  The solutions
which do not have this symmetry are denoted with {\em asym} in
the table.

\begin{table}[htb]
{\begin{center}
\begin{tabular}{| c | c |  c | c | c |}
    \hline 
   Degree of cardinal  & Number of   & Degree of Exact   &         &  \\
   functions (d)            & Points (N)  & Integration (d+e) &  Error  & Notes \\ \hline
       1  &  3 &  2 & 4.4 $\times 10^{-16}$ &  \\ 
       2  &  6 &  4 & 9.7 $\times 10^{-16}$ &  \\ 
       3  & 10 &  5 & 1.7 $\times 10^{-14}$ &  \\ 
       4  & 15 &  7 & 2.1 $\times 10^{-14}$ &  \\ 
       5  & 21 &  9 & 2.8 $\times 10^{-14}$ &  \\ 
       6  & 28 & 11 & 4.7 $\times 10^{-15}$ & asym  \\ 
       7  & 36 & 13 & 2.2 $\times 10^{-14}$ & asym,new \\ 
       8  & 45 & 14 & 1.8 $\times 10^{-15}$ &  \\ 
       9  & 55 & 16 & 8.6 $\times 10^{-15}$ & asym,new  \\ 
      10  & 66 & 18 & 3.3 $\times 10^{-14}$ & asym,new  \\ 
      11  & 78 & 20 & 2.8 $\times 10^{-14}$ & asym,new  \\ 
      12  & 91 & 21 & 2.9 $\times 10^{-14}$ & new \\ 
      13  & 105& 23 & 3.3 $\times 10^{-14}$ & new \\ 
      14  & 120& 25 & 4.3 $\times 10^{-14}$ & asym,new \\ 
\hline
\end{tabular}
\caption{Quadrature points computed with the cardinal function algorithm.  
In all cases, the quadrature weights are positive and the points are not
outside the triangle.  Solutions which are not $D_3$ symmetric are
denoted by {\em asym}.  Solutions which improve upon previously
published results are denoted by {\em new}.
}
\end{center}}
\label{T:quad}
\end{table}

Quadrature formulas denoted by {\em new} in the table represent 
formulas which improve upon the best previously published results,
as taken from the extensive database described in \cite{Cools2003} as
well as the quadrature points presented in \cite{WandzuraXiao2003}
(which are not included in the database as of this writing).  
The new solutions for integration degree $d+e$ from 16 to 25 have fewer
points then the previously published results.  
For $d+e=13$, the previous result with the fewest number 
of quadrature points \cite{BerntsenEspelid1990} has $N=36$.  
The points in \cite{BerntsenEspelid1990} are symmetric, but some are outside
the triangle and not all weights are positive.  
The result presented here also has 36 points, all of
which are inside the triangle, the weights are positive, but
the points are asymmetric.

Plots for all of the quadrature points are given in the
figures.  For the plots, the right triangle has been
mapped lineally to the equilateral triangle in order to make
the asymmetry in the points more visible.  
The coordinates of the quadrature points are given in
Appendex A.  They are available electronically by downloading
the \TeX\  source of this paper from the arXiv.

\section{Summary}
\label{section:summary}

We have presented results from a cardinal function algorithm for computing
multi-variate quadrature points.   The algorithm was applied
to the triangle, where optimal (in the sense of Equations~\ref{E:optimal1}
and~\ref{E:optimal2}) formulas
were constructed for integrating polynomials up to degree 25.  
Seven of these quadrature
formulas improve on previously known results.

\bibliography{../../../tex/genbib,../../../tex/Quadrature}

\appendix
\clearpage
\section{Plots of quadrature points}

\begin{figure}[!h]
\begin{center}
\includegraphics[width=2.2in]{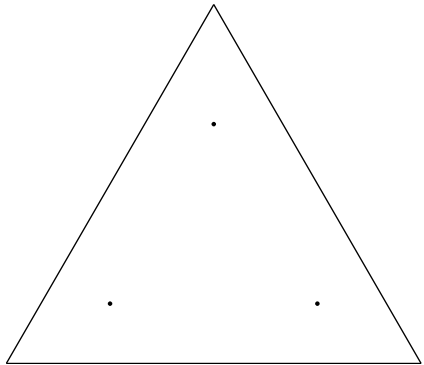}
\caption{Quadrature points for the triangle which
exactly integrate polynomials of degree 2.
}
\end{center}
\end{figure}

\begin{figure}[!h]
\begin{center}
\includegraphics[width=2.2in]{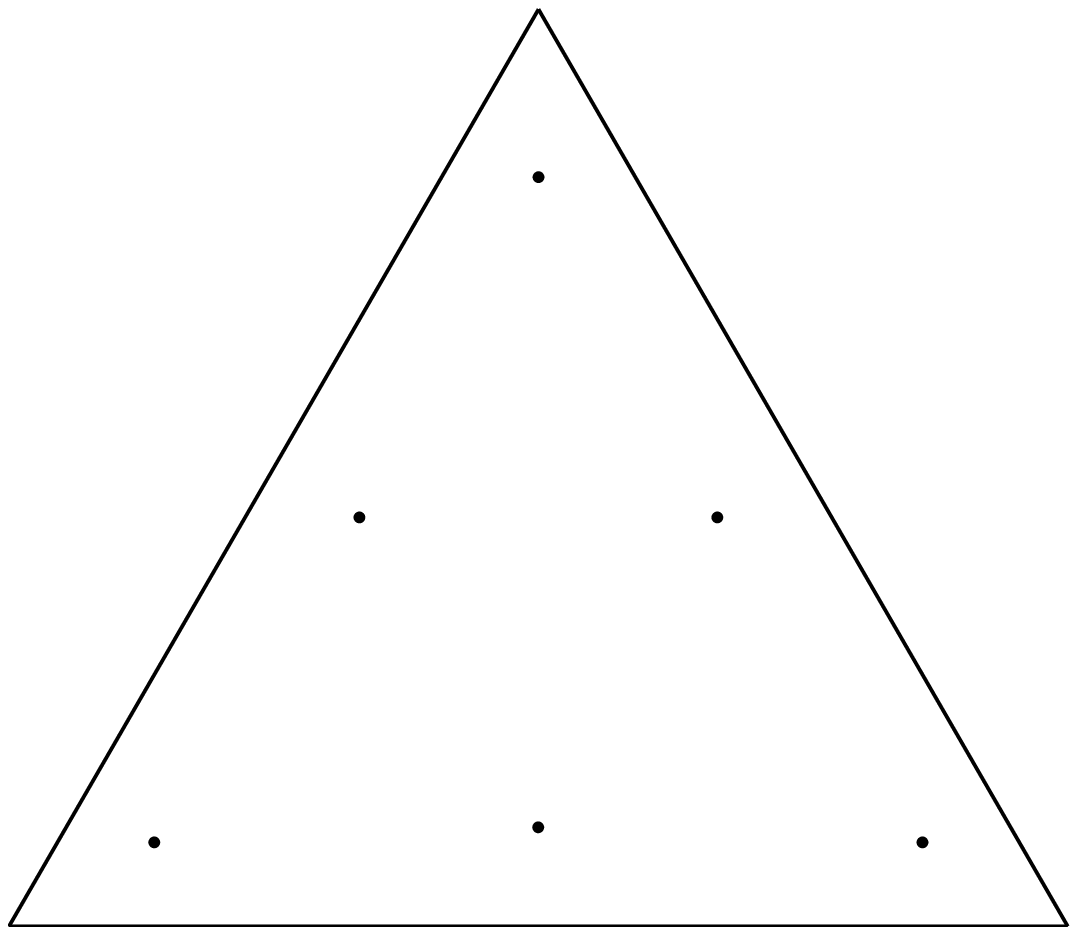}
\caption{Quadrature points for the triangle which
exactly integrate polynomials of degree 4.
}
\end{center}
\end{figure}

\begin{figure}[!h]
\begin{center}
\includegraphics[width=2.2in]{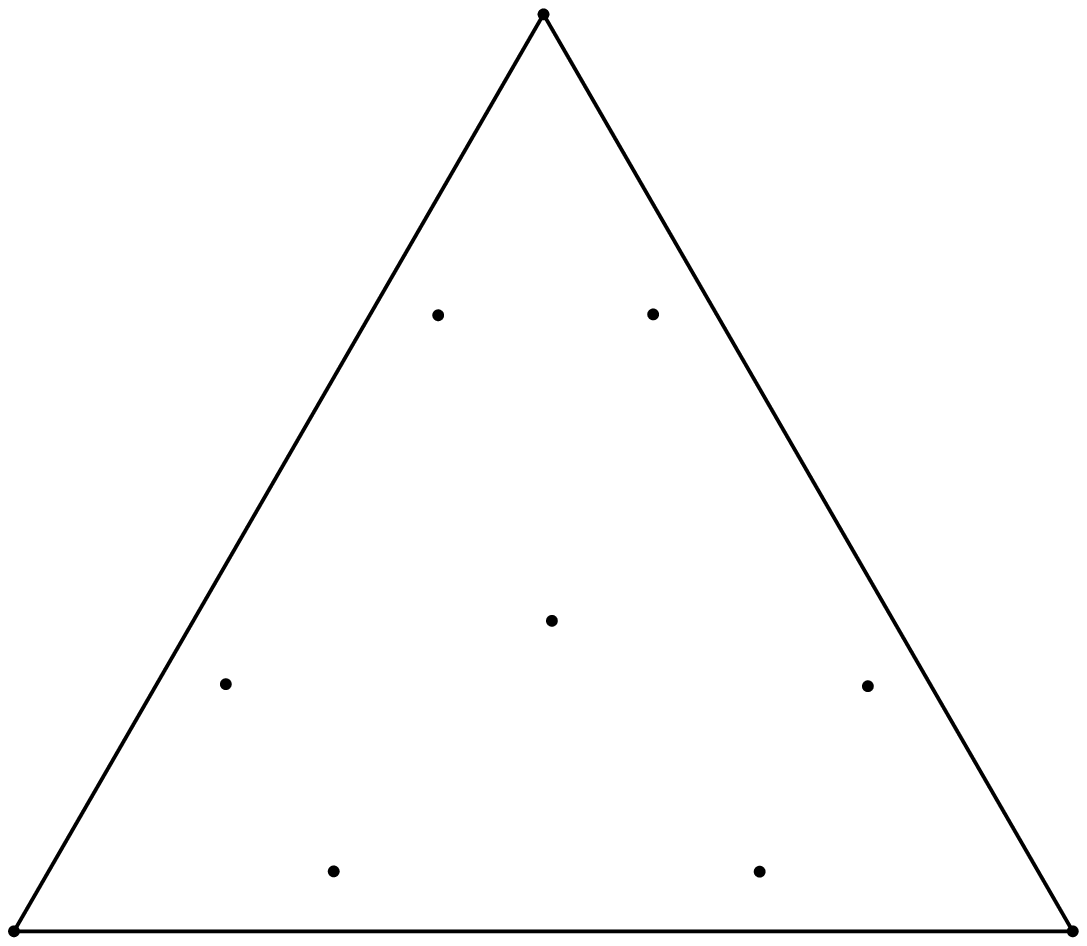}
\caption{Quadrature points for the triangle which
exactly integrate polynomials of degree 5.
}
\end{center}
\end{figure}

\begin{figure}[!h]
\begin{center}
\includegraphics[width=2.2in]{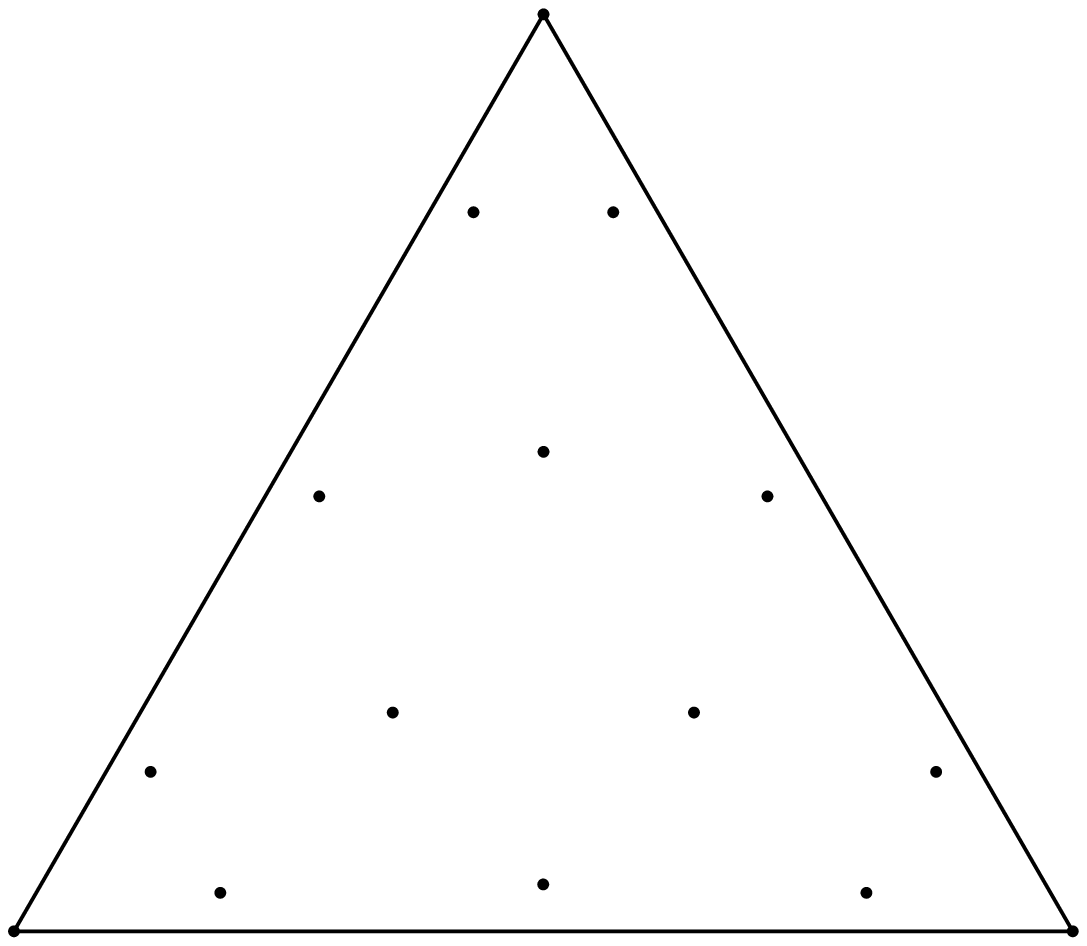}
\caption{Quadrature points for the triangle which
exactly integrate polynomials of degree 7.
}
\end{center}
\end{figure}

\begin{figure}[!h]
\begin{center}
\includegraphics[width=2.2in]{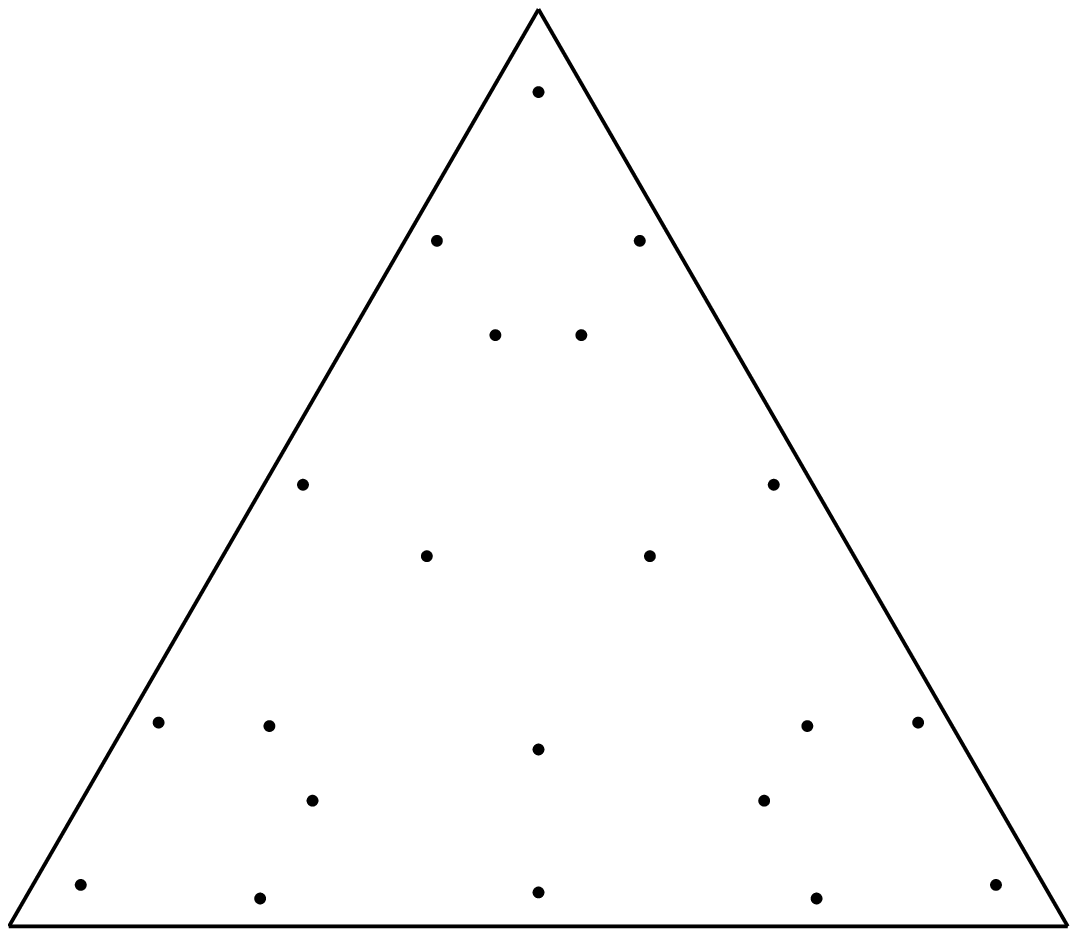}
\caption{Quadrature points for the triangle which
exactly integrate polynomials of degree 9.
}
\end{center}
\end{figure}

\begin{figure}[!h]
\begin{center}
\includegraphics[width=2.2in]{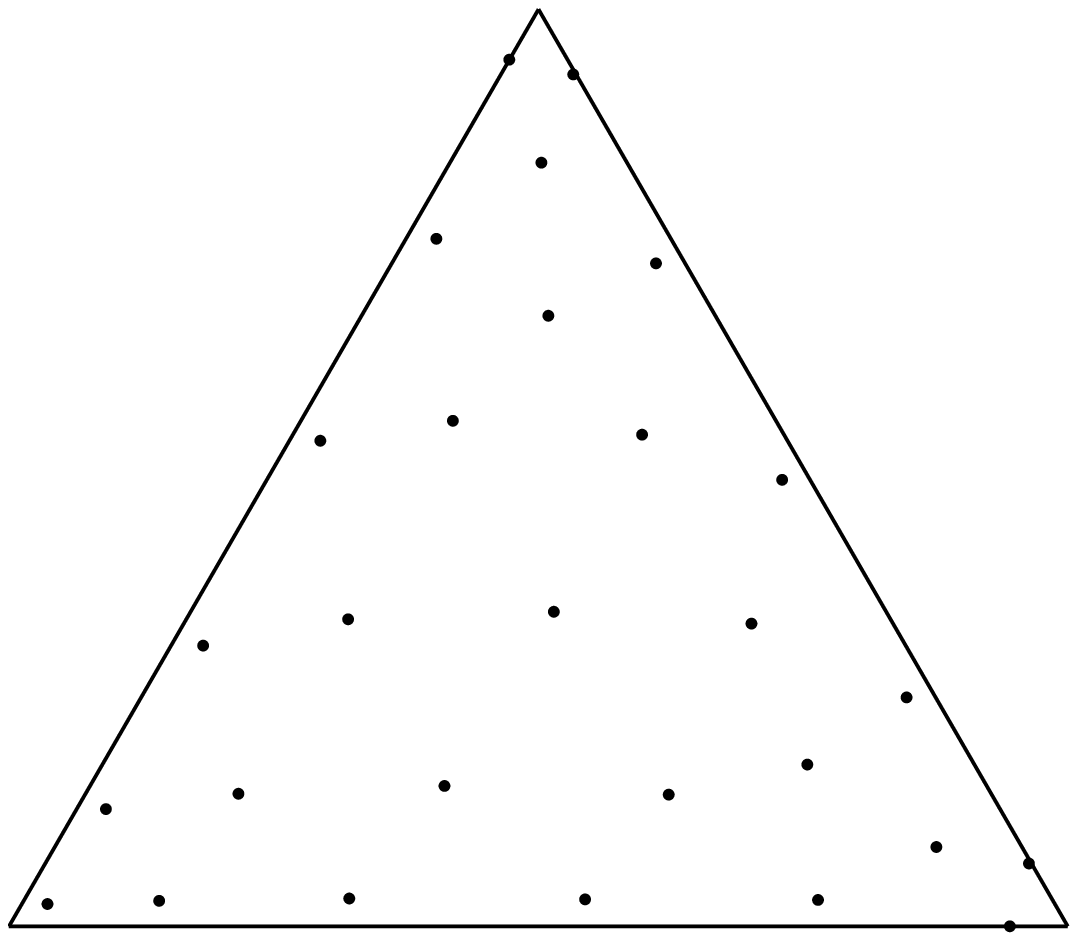}
\caption{Quadrature points for the triangle which
exactly integrate polynomials of degree 11.
}
\end{center}
\end{figure}

\begin{figure}[!h]
\begin{center}
\includegraphics[width=2.2in]{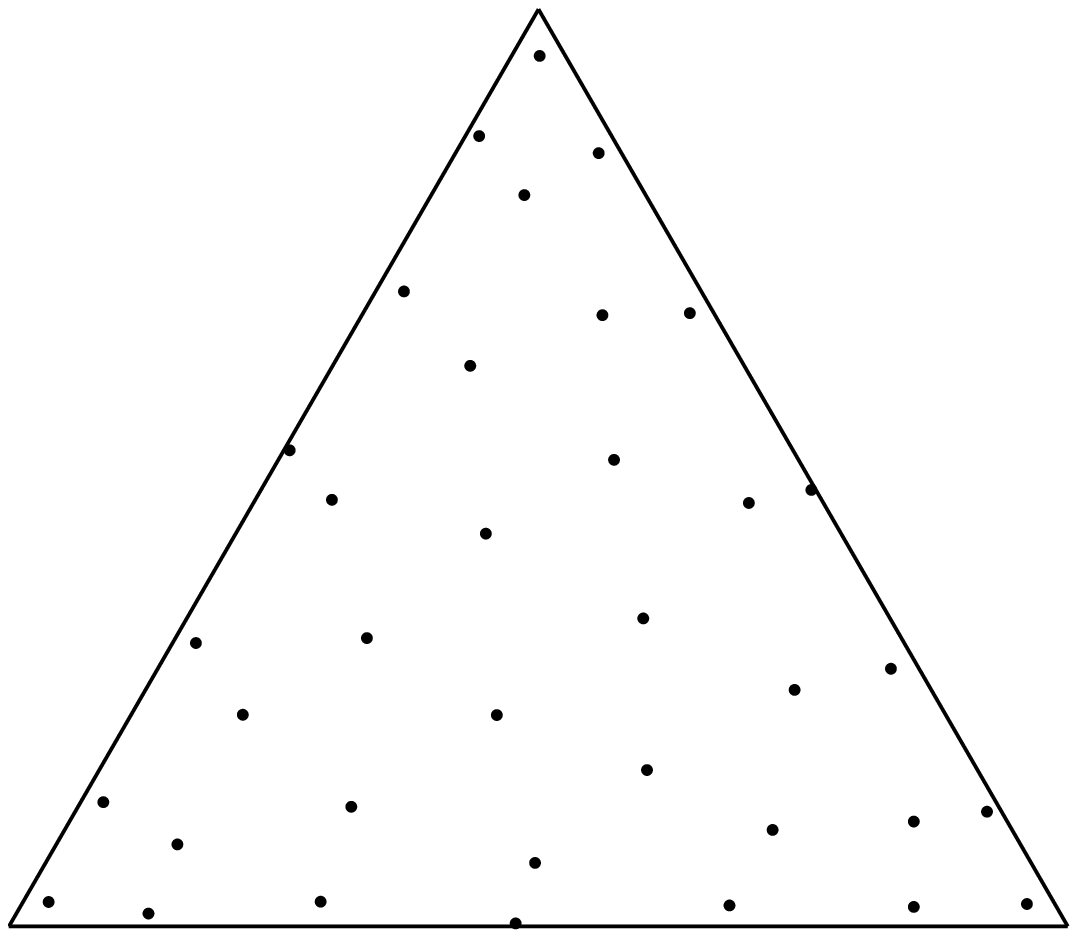}
\caption{Quadrature points for the triangle which
exactly integrate polynomials of degree 13.
}
\end{center}

\end{figure}
\begin{figure}[!h]
\begin{center}
\includegraphics[width=2.2in]{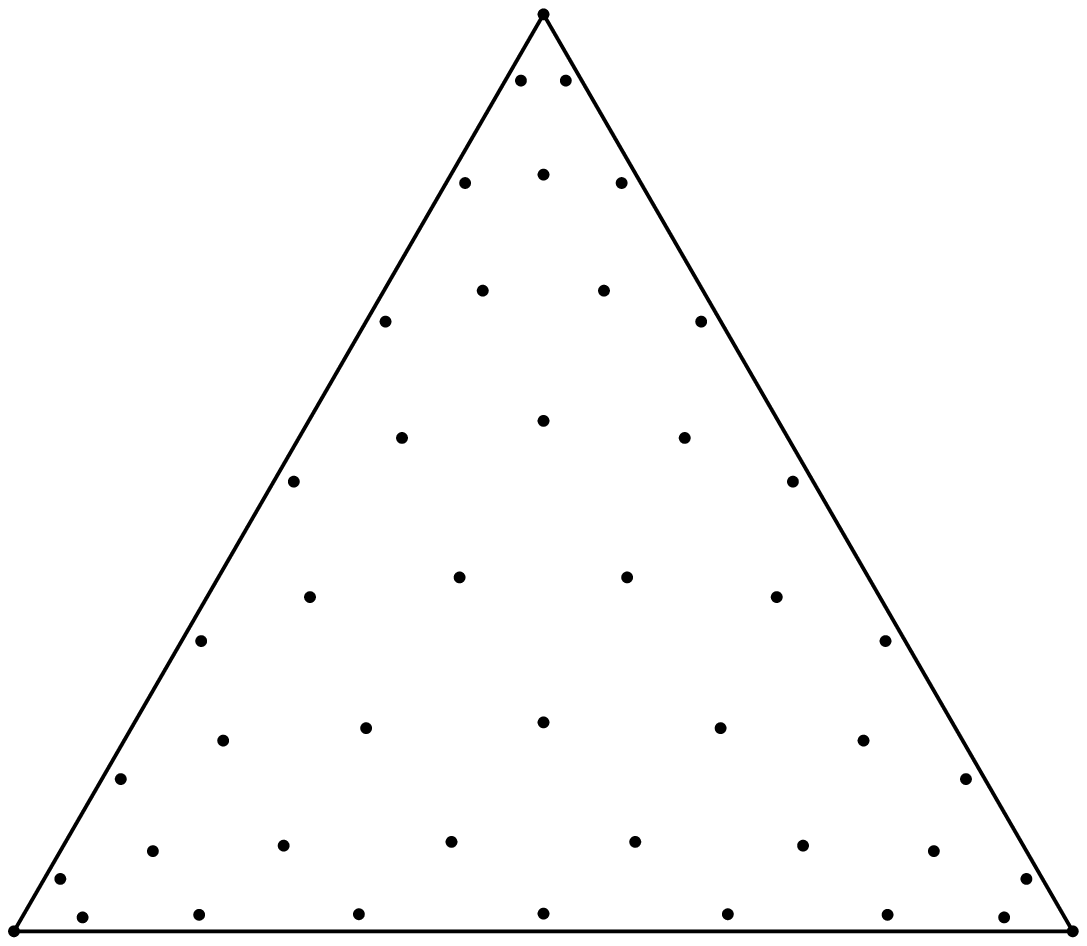}
\caption{Quadrature points for the triangle which
exactly integrate polynomials of degree 14.
}
\end{center}
\end{figure}

\begin{figure}[!h]
\begin{center}
\includegraphics[width=2.2in]{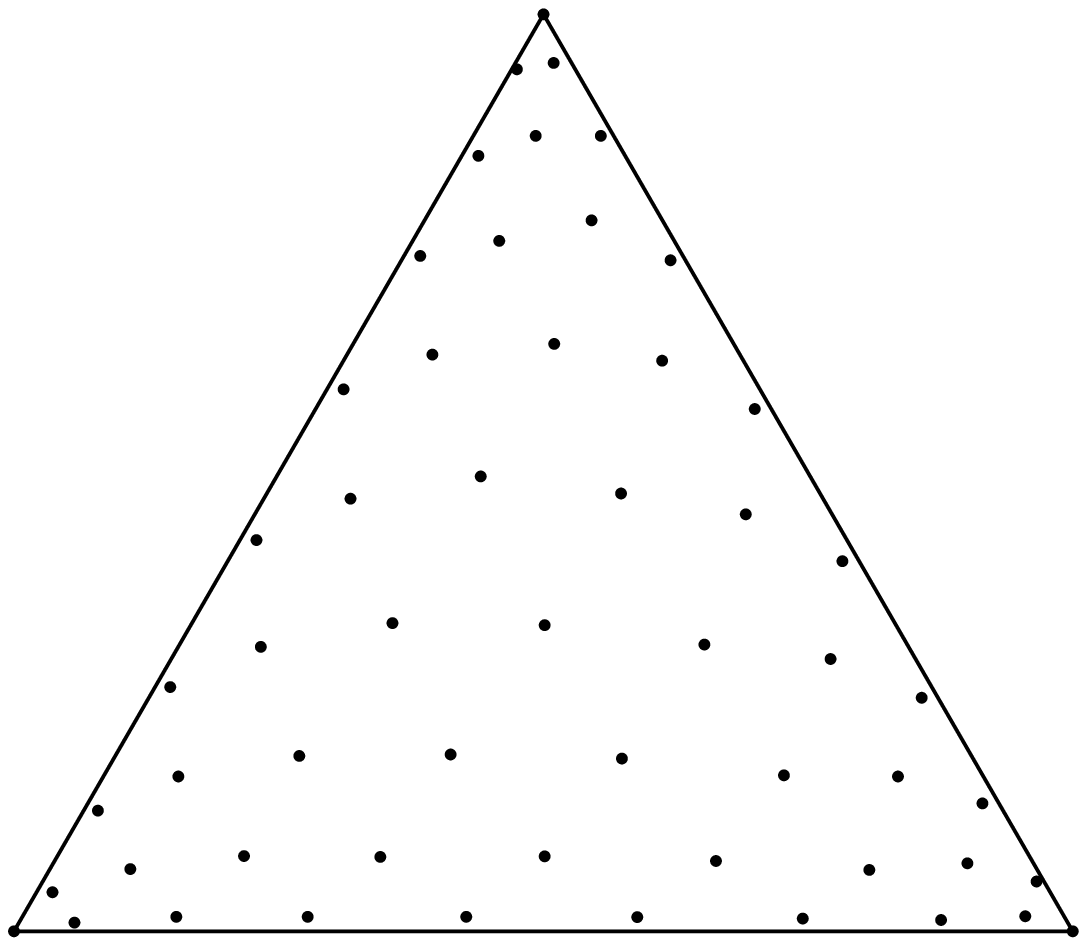}
\caption{Quadrature points for the triangle which
exactly integrate polynomials of degree 16.
}
\end{center}
\end{figure}

\begin{figure}[!h]
\begin{center}
\includegraphics[width=2.2in]{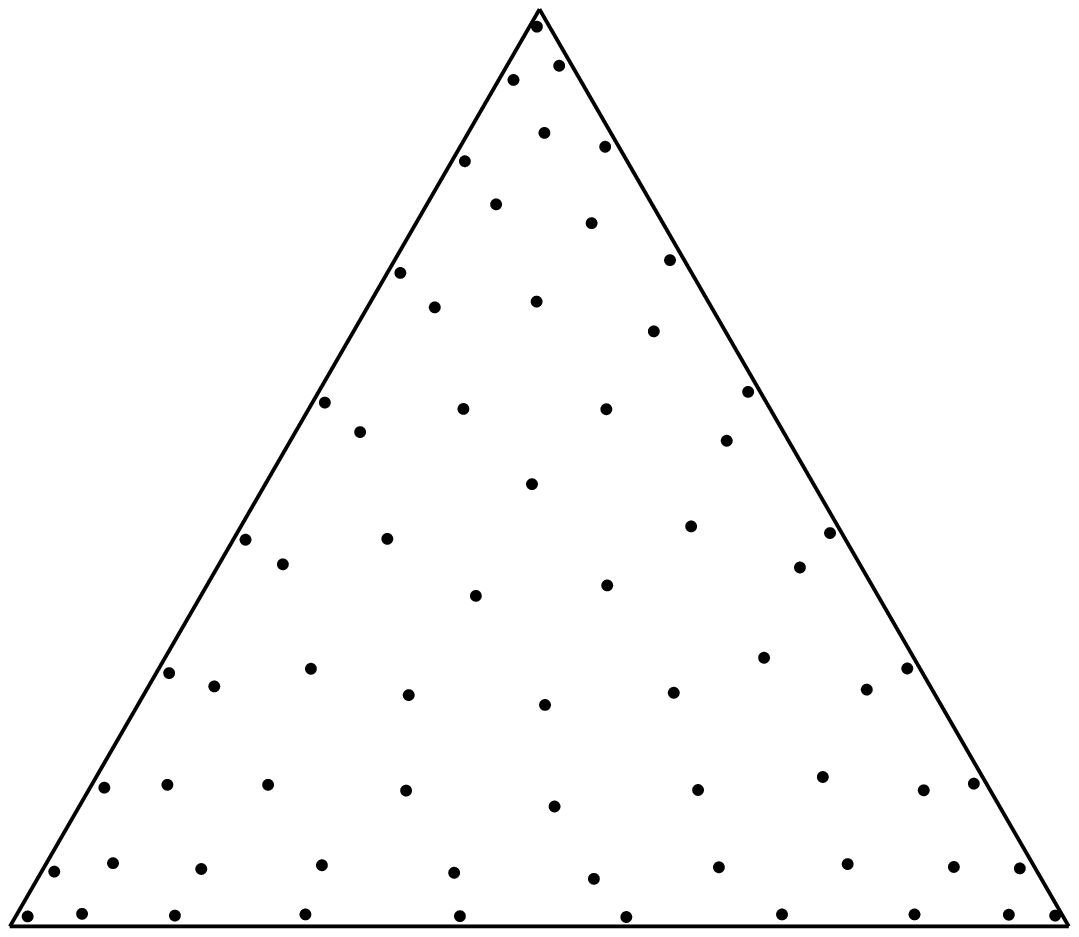}
\caption{Quadrature points for the triangle which
exactly integrate polynomials of degree 18.
}
\end{center}
\end{figure}

\begin{figure}[!h]
\begin{center}
\includegraphics[width=2.2in]{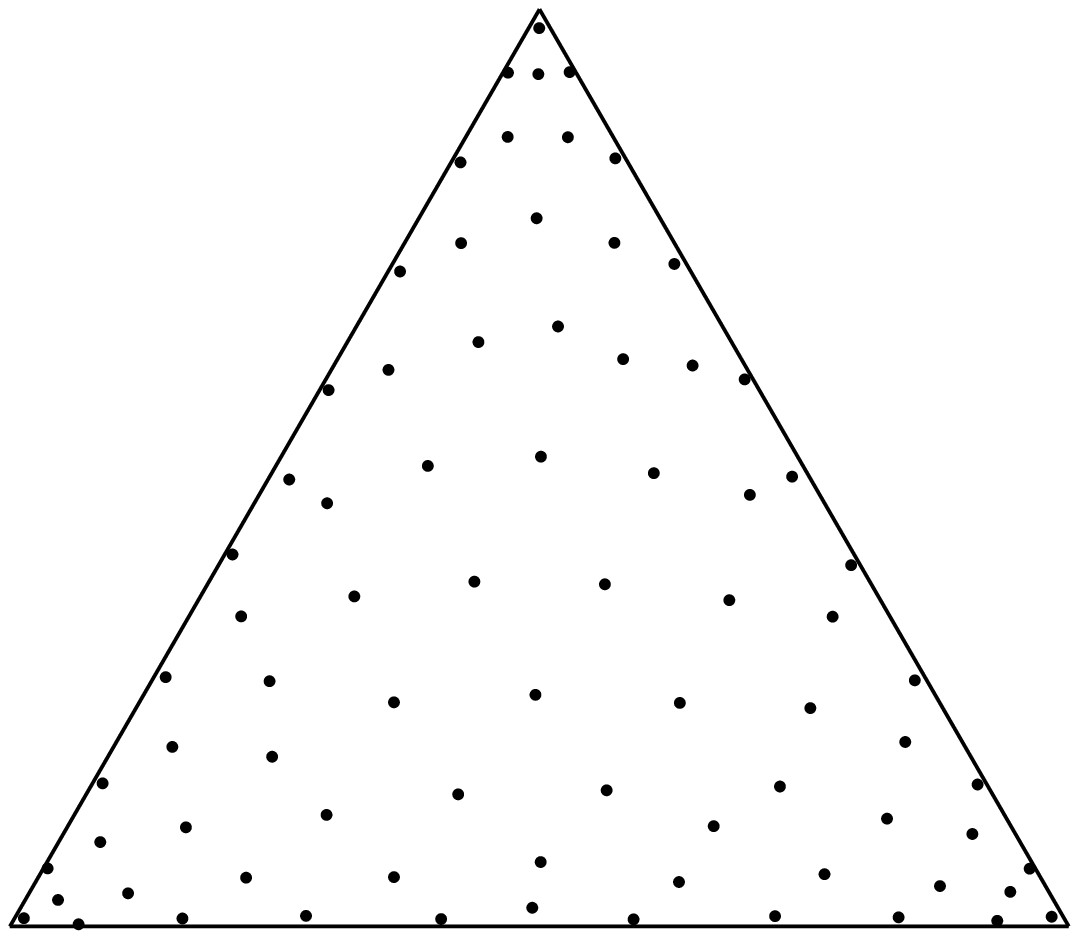}
\caption{Quadrature points for the triangle which
exactly integrate polynomials of degree 20.
}
\end{center}
\end{figure}

\begin{figure}[!h]
\begin{center}
\includegraphics[width=2.2in]{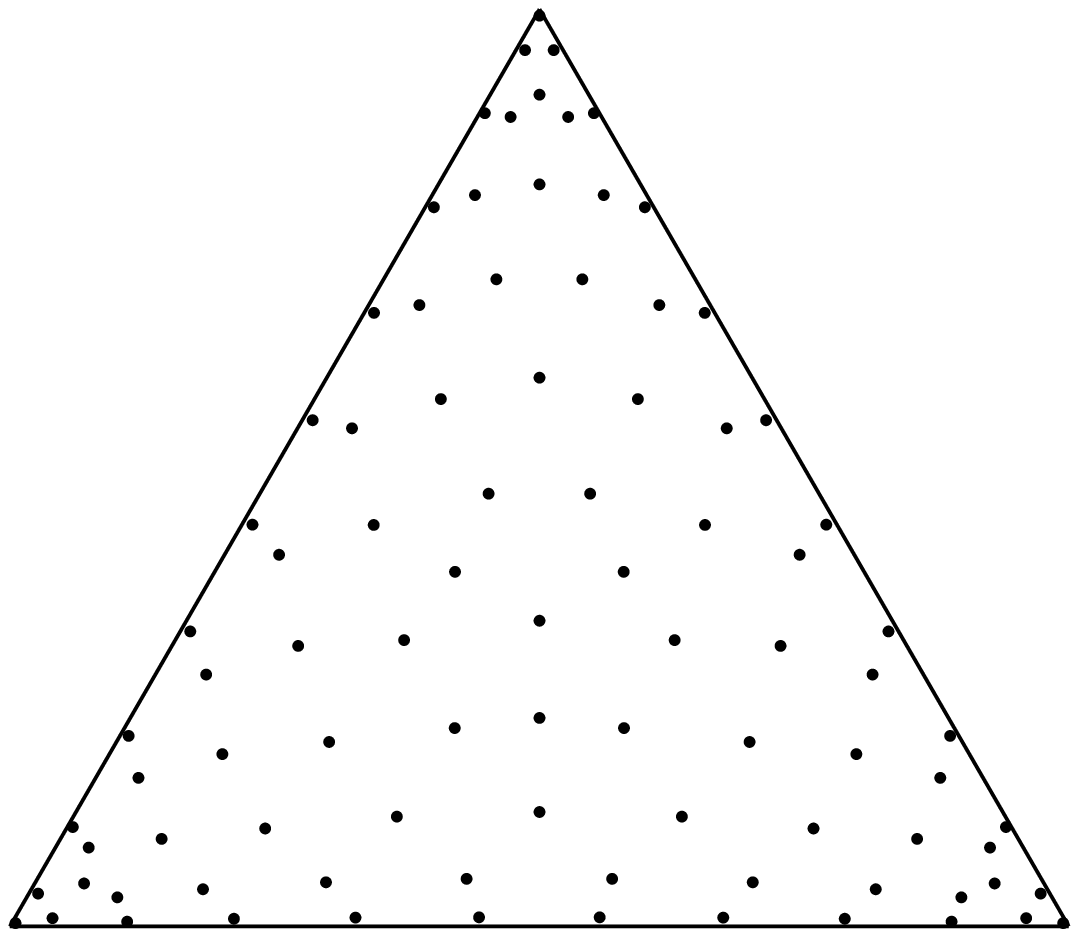}
\caption{Quadrature points for the triangle which
exactly integrate polynomials of degree 21.
}
\end{center}
\end{figure}

\begin{figure}[!h]
\begin{center}
\includegraphics[width=2.2in]{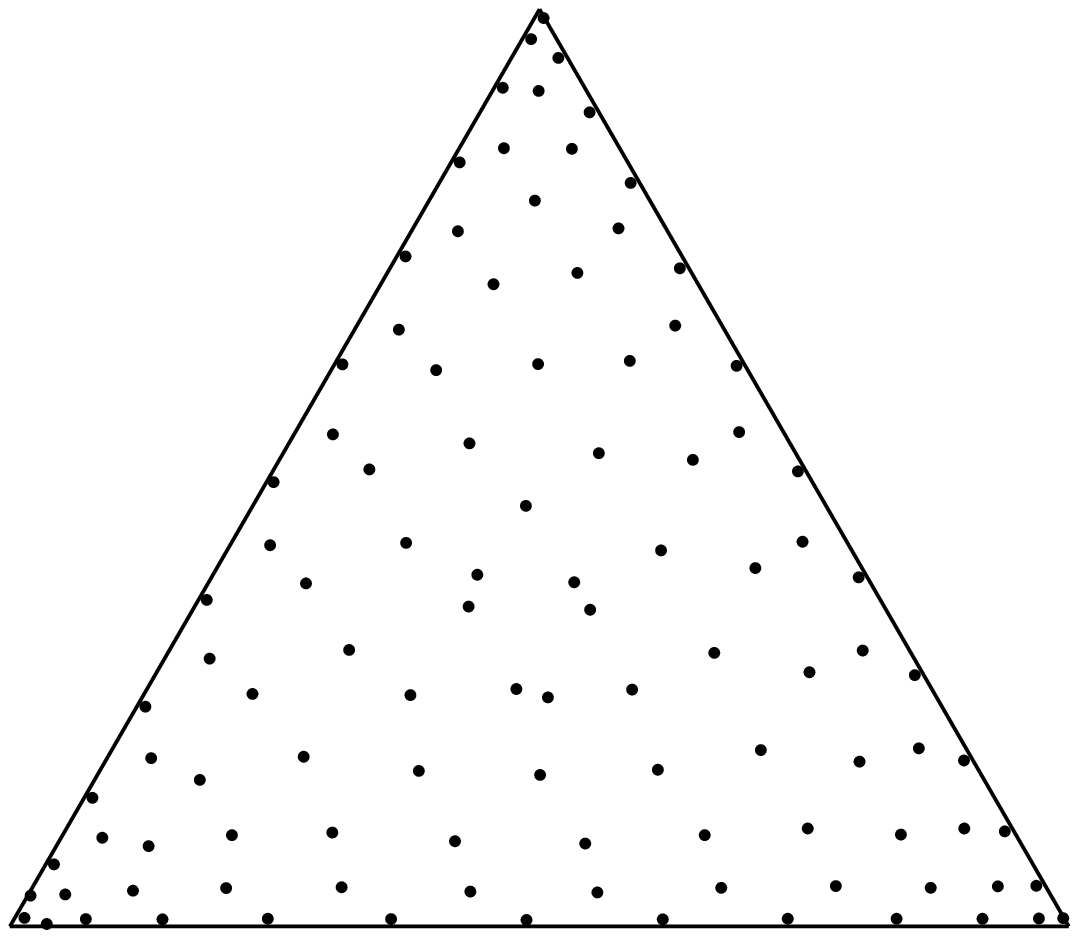}
\caption{Quadrature points for the triangle which
exactly integrate polynomials of degree 23.}
\end{center}
\end{figure}

\begin{figure}[!h]
\begin{center}
\includegraphics[width=2.2in]{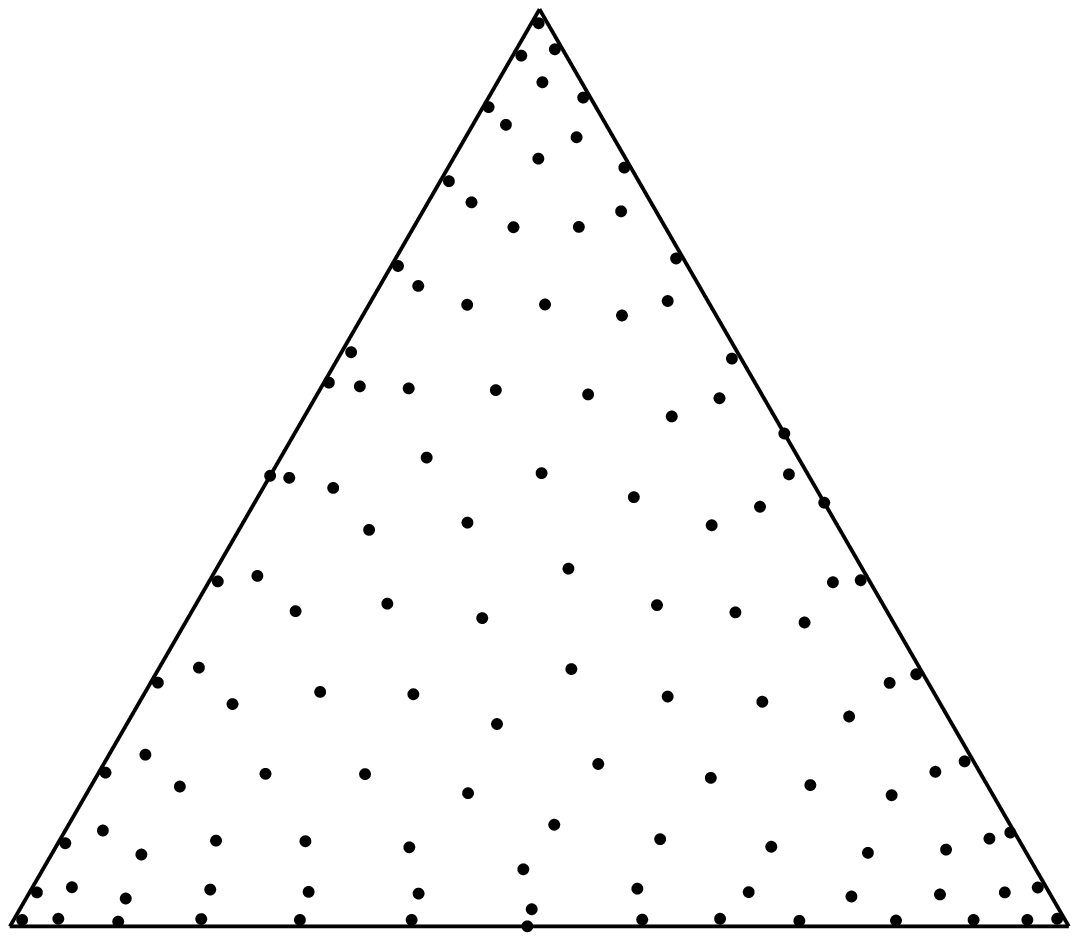}
\caption{Quadrature points for the triangle which
exactly integrate polynomials of degree 25.}
\end{center}
\end{figure}

\clearpage
\section{Tables of quadrature points}
We now list the coordinates of the quadrature points described in
Table~\ref{T:quad}. 
For each line, we give the first two barycentric coordinates of each point
(equivalent to the $x$ and $y$ coordinates after an equilateral
triangle is linearly mapped to the unit right triangle 
$x \ge 0$,
$y\ge 0$ and $x+y \le 1$) followed by the associated quadrature weight. 
The third
barycentric coordinate is defined such that the sum of all three
coordinates is one.   These points are available electronically by downloading
the \TeX\  source of this paper from the arXiv.  

\begin{multicols}{2}{\centerline{Tables of Points}}
\begin{ttfamily}
\scriptsize
\obeylines
\parindent 0pt

\input coords.txt

\end{ttfamily}
\end{multicols}

\end{document}

%% file: coords.txt
integration degree=2  N=3:
 0.1666666666667   0.6666666666667  0.6666666666667
 0.6666666666667   0.1666666666667  0.6666666666667
 0.1666666666667   0.1666666666667  0.6666666666667
integration degree=4  N=6:
 0.0915762135098   0.0915762135098  0.2199034873106
 0.8168475729805   0.0915762135098  0.2199034873106
 0.0915762135098   0.8168475729805  0.2199034873106
 0.1081030181681   0.4459484909160  0.4467631793560
 0.4459484909160   0.1081030181681  0.4467631793560
 0.4459484909160   0.4459484909160  0.4467631793560
integration degree=5  N=10:
 0.0000000000000   1.0000000000000  0.0262712099504
 1.0000000000000   0.0000000000000  0.0262716612068
 0.0000000000000   0.0000000000000  0.0274163947600
 0.2673273531185   0.6728199218710  0.2348383865823
 0.6728175529461   0.2673288599482  0.2348412238268
 0.0649236350054   0.6716530111494  0.2480251793114
 0.6716498539042   0.0649251690029  0.2480304922521
 0.0654032456800   0.2693789366453  0.2518604605529
 0.2693767069140   0.0654054874919  0.2518660533658
 0.3386738503896   0.3386799893027  0.4505789381914
integration degree=7  N=15:
 1.0000000000000   0.0000000000000  0.0102558174092
 0.0000000000000   0.0000000000000  0.0102558174092
 0.0000000000000   1.0000000000000  0.0102558174092
 0.7839656651012   0.0421382841642  0.1116047046647
 0.1738960507345   0.7839656651012  0.1116047046647
 0.1738960507345   0.0421382841642  0.1116047046647
 0.0421382841642   0.1738960507345  0.1116047046647
 0.7839656651012   0.1738960507345  0.1116047046647
 0.0421382841642   0.7839656651012  0.1116047046647
 0.4743880861752   0.4743880861752  0.1679775595335
 0.4743880861752   0.0512238276497  0.1679775595335
 0.0512238276497   0.4743880861752  0.1679775595335
 0.2385615300181   0.5228769399639  0.2652238803946
 0.5228769399639   0.2385615300181  0.2652238803946
 0.2385615300181   0.2385615300181  0.2652238803946
integration degree=9  N=21:
 0.0451890097844   0.0451890097844  0.0519871420646
 0.0451890097844   0.9096219804312  0.0519871420646
 0.9096219804312   0.0451890097844  0.0519871420646
 0.7475124727339   0.0304243617288  0.0707034101784
 0.2220631655373   0.0304243617288  0.0707034101784
 0.7475124727339   0.2220631655373  0.0707034101784
 0.2220631655373   0.7475124727339  0.0707034101784
 0.0304243617288   0.7475124727339  0.0707034101784
 0.0304243617288   0.2220631655373  0.0707034101784
 0.1369912012649   0.2182900709714  0.0909390760952
 0.6447187277637   0.2182900709714  0.0909390760952
 0.1369912012649   0.6447187277637  0.0909390760952
 0.2182900709714   0.6447187277637  0.0909390760952
 0.2182900709714   0.1369912012649  0.0909390760952
 0.6447187277637   0.1369912012649  0.0909390760952
 0.0369603304334   0.4815198347833  0.1032344051380
 0.4815198347833   0.0369603304334  0.1032344051380
 0.4815198347833   0.4815198347833  0.1032344051380
 0.4036039798179   0.1927920403641  0.1881601469167
 0.4036039798179   0.4036039798179  0.1881601469167
 0.1927920403641   0.4036039798179  0.1881601469167
integration degree=11  N=28:
 0.0000000000000   0.9451704450174  0.0114082494033
 0.9451704450173   0.0000000000000  0.0114082494033
 0.9289002405719   0.0685505797224  0.0132691285720
 0.0685505797224   0.9289002405717  0.0132691285720
 0.0243268355615   0.0243268355616  0.0155865773350
 0.1279662835335   0.0277838749488  0.0408274780428
 0.0277838749488   0.1279662835337  0.0408274780429
 0.0287083428360   0.7498347588657  0.0579849665116
 0.7498347588656   0.0287083428360  0.0579849665116
 0.7228007909707   0.2497602062385  0.0601385247663
 0.2497602062386   0.7228007909707  0.0601385247663
 0.0865562992839   0.8325513856997  0.0625273888433
 0.8325513856998   0.0865562992839  0.0625273888433
 0.3061619157672   0.0303526617491  0.0639684321504
 0.0303526617491   0.3061619157675  0.0639684321504
 0.4868610595047   0.4868610595047  0.0661325872161
 0.6657904293017   0.1765456154219  0.0668503236820
 0.1765456154221   0.6657904293016  0.0668503236821
 0.0293121007360   0.5295657488669  0.0686904305977
 0.5295657488667   0.0293121007360  0.0686904305977
 0.1444673824391   0.1444673824391  0.1002717543859
 0.3299740111411   0.5361815729050  0.1143136784099
 0.5361815729052   0.3299740111409  0.1143136784099
 0.5511507516862   0.1437790861923  0.1223648146752
 0.1437790861923   0.5511507516862  0.1223648146752
 0.3348066587327   0.1529619437161  0.1394422334178
 0.1529619437161   0.3348066587327  0.1394422334178
 0.3430183498147   0.3430183498147  0.1744377829182
integration degree=13  N=36:
 0.0242935351590   0.9493059293846  0.0166240998757
 0.0265193427722   0.0242695130640  0.0166811699778
 0.9492126023551   0.0265067966437  0.0166830569067
 0.0033775763749   0.4767316412363  0.0175680870083
 0.4757672298101   0.5198921829102  0.0184474661845
 0.5190783193471   0.0055912706202  0.0197942410188
 0.8616839745321   0.0133996048618  0.0203540395855
 0.1249209759926   0.8613054321334  0.0206852863940
 0.0138565453861   0.1247733717358  0.0208271366086
 0.0211887064222   0.8438438351223  0.0317819778279
 0.8432296787219   0.1354563645830  0.0320472035241
 0.1354231797865   0.0213482820656  0.0320607681146
 0.3088853510679   0.0221919663014  0.0430765959183
 0.6685057595169   0.3089012879389  0.0438473415339
 0.0226545012557   0.6691709943321  0.0439209672733
 0.2808515408772   0.6924718155106  0.0479951923691
 0.6922446749051   0.0268723345026  0.0483806260733
 0.0268617447119   0.2810093973222  0.0484867423375
 0.1141778485470   0.7973581413586  0.0556964488024
 0.7974807922061   0.0879806508791  0.0561026364356
 0.0892807293894   0.1145020561128  0.0565190123693
 0.1052487892455   0.6686904119922  0.0689289890670
 0.6663022280740   0.2275051631832  0.0717213336089
 0.2307803737547   0.1054572561221  0.0727453920976
 0.1705059157540   0.5174064398658  0.0788807336737
 0.5086593973043   0.3170523855209  0.0810114345512
 0.3141823862281   0.1810706361659  0.0825725299055
 0.4617460817864   0.4678594539804  0.0842044567330
 0.0693087496081   0.4622856042085  0.0843585533305
 0.4651955259268   0.0724357805669  0.0851969868488
 0.2578625857893   0.6131395039177  0.0902845328052
 0.6112627766779   0.1300360834609  0.0914283143485
 0.1305182135934   0.2581713828884  0.0916279065409
 0.4281437991828   0.2362005969817  0.1025573374896
 0.3356995783730   0.4311026308588  0.1033159661413
 0.2305424298836   0.3456013949376  0.1035854367193
integration degree=14  N=45:
 0.0000000000000   1.0000000000000  0.0010616711990
 1.0000000000000   0.0000000000000  0.0010616711990
 0.0000000000000   0.0000000000000  0.0010616711990
 0.0573330873026   0.0151382269814  0.0131460236101
 0.0573330873026   0.9275286857160  0.0131460236101
 0.9275286857160   0.0573330873026  0.0131460236101
 0.0151382269814   0.0573330873026  0.0131460236101
 0.9275286857160   0.0151382269814  0.0131460236101
 0.0151382269814   0.9275286857160  0.0131460236101
 0.8159625040711   0.1659719969565  0.0242881926949
 0.8159625040711   0.0180654989724  0.0242881926949
 0.1659719969565   0.8159625040711  0.0242881926949
 0.0180654989724   0.8159625040711  0.0242881926949
 0.1659719969565   0.0180654989724  0.0242881926949
 0.0180654989724   0.1659719969565  0.0242881926949
 0.3165475556378   0.0186886898773  0.0316799866332
 0.6647637544849   0.0186886898773  0.0316799866332
 0.0186886898773   0.6647637544849  0.0316799866332
 0.0186886898773   0.3165475556378  0.0316799866332
 0.3165475556378   0.6647637544849  0.0316799866332
 0.6647637544849   0.3165475556378  0.0316799866332
 0.0192662192492   0.4903668903754  0.0349317947036
 0.4903668903754   0.0192662192492  0.0349317947036
 0.4903668903754   0.4903668903754  0.0349317947036
 0.0875134669581   0.8249730660837  0.0383664533945
 0.0875134669581   0.0875134669581  0.0383664533945
 0.8249730660837   0.0875134669581  0.0383664533945
 0.0935526036219   0.2079865423167  0.0578369491210
 0.0935526036219   0.6984608540613  0.0578369491210
 0.2079865423167   0.0935526036219  0.0578369491210
 0.6984608540613   0.0935526036219  0.0578369491210
 0.6984608540613   0.2079865423167  0.0578369491210
 0.2079865423167   0.6984608540613  0.0578369491210
 0.0974892983467   0.5380088595149  0.0725821687394
 0.3645018421383   0.0974892983467  0.0725821687394
 0.5380088595149   0.0974892983467  0.0725821687394
 0.5380088595149   0.3645018421383  0.0725821687394
 0.3645018421383   0.5380088595149  0.0725821687394
 0.0974892983467   0.3645018421383  0.0725821687394
 0.2217145894873   0.5565708210253  0.0897856524107
 0.5565708210253   0.2217145894873  0.0897856524107
 0.2217145894873   0.2217145894873  0.0897856524107
 0.3860471669296   0.2279056661408  0.1034544533617
 0.2279056661408   0.3860471669296  0.1034544533617
 0.3860471669296   0.3860471669296  0.1034544533617
integration degree=16  N=55:
 1.0000000000000   0.0000000000000  0.0006202599851
 0.0000000000000   1.0000000000000  0.0006315174712
 0.0000000000000   0.0000000000000  0.0007086601559
 0.9398863583577   0.0049848744634  0.0055163716168
 0.0543806683058   0.9386405618617  0.0062692407656
 0.0093940049164   0.0526424462697  0.0078531408826
 0.0164345086362   0.9469035517351  0.0094551483864
 0.9469487269862   0.0363373677167  0.0097824511271
 0.0426604005768   0.0151224541799  0.0099861643489
 0.0122269495439   0.8693773510664  0.0137553818816
 0.8673696521047   0.1204917285774  0.0140979178040
 0.8456744021389   0.0157763967870  0.0149646864337
 0.1395759632103   0.8448120870375  0.0156097503612
 0.1317821743231   0.0135009605584  0.0157683693348
 0.0157955126300   0.1455274938536  0.0175794546383
 0.7365462884436   0.0155697540908  0.0204113840270
 0.0139688430330   0.7379836894450  0.0209562878616
 0.2547895186039   0.7297615689771  0.0210713412998
 0.7316386522555   0.2543076683315  0.0217646760202
 0.0157253728951   0.2696239795791  0.0222288408699
 0.2662302843647   0.0144783956308  0.0224186693682
 0.8673504065214   0.0591679410400  0.0230122616993
 0.0741493666957   0.8634782575061  0.0236813902500
 0.0159285948360   0.4191238955238  0.0257464643368
 0.0156061028068   0.5809222921146  0.0257956801608
 0.5910094817484   0.0159251452651  0.0258072327610
 0.4034771496889   0.5806700368104  0.0260343232059
 0.5694745628526   0.4149495146302  0.0265768141609
 0.0678493700650   0.0761218678591  0.0265784761831
 0.4265968590272   0.0157509692312  0.0267532329238
 0.0670982507890   0.7741898312421  0.0375787806641
 0.7528310231480   0.0819119495639  0.0383065894195
 0.7753727783557   0.1577128457292  0.0384849695025
 0.1689073157787   0.7503943099742  0.0389619825852
 0.1687335832919   0.0708311507268  0.0394604111547
 0.0821244708436   0.1762996626771  0.0412364778098
 0.6288705363345   0.0807744953317  0.0512872438483
 0.0811413015266   0.3054373589776  0.0516405641935
 0.2969112065080   0.6227485988871  0.0518230042269
 0.0767542314171   0.6247247149546  0.0528527988181
 0.6223022333845   0.3011485821166  0.0538505573027
 0.3103786288051   0.0779098365079  0.0541895329319
 0.0819218215187   0.4603633038351  0.0584737146444
 0.4717022665013   0.0821554006797  0.0592863168363
 0.4546603415250   0.4637565033890  0.0594358276749
 0.1701091339237   0.6422277808188  0.0631800255863
 0.6406004329487   0.1898293537256  0.0632926845153
 0.1912267583717   0.1739955685343  0.0640707361772
 0.1885315767070   0.4798914070406  0.0812040595918
 0.4772929957691   0.3348356598119  0.0814437513530
 0.3126974621760   0.4957972197259  0.0814679201241
 0.4961225945946   0.1927553668904  0.0815050548084
 0.1928805312867   0.3161015807261  0.0815164664939
 0.3360041453816   0.1894892801290  0.0816931059623
 0.3337280550848   0.3343571021811  0.0923218334531
integration degree=18  N=66:
 0.0116731059668   0.9812565951289  0.0025165756986
 0.9810030858388   0.0071462504863  0.0025273452007
 0.0106966317092   0.0115153933376  0.0033269295333
 0.9382476983551   0.0495570591341  0.0081503492125
 0.0126627518417   0.9370123620615  0.0086135525742
 0.0598109409984   0.0121364578922  0.0087786746179
 0.0137363297927   0.0612783625597  0.0097099585562
 0.9229527959405   0.0141128270602  0.0102466211915
 0.0633107354993   0.9220197291727  0.0108397688341
 0.0117265100335   0.1500520475229  0.0129385390176
 0.1554720587323   0.8325147121589  0.0136339823583
 0.8343293888982   0.0125228158759  0.0138477328147
 0.8501638031957   0.1371997508736  0.0139421540105
 0.0128816350522   0.8477627063479  0.0144121399968
 0.1510801608959   0.0136526924039  0.0153703455534
 0.0101917879217   0.5770438618345  0.0162489802253
 0.2813372399303   0.7066853759623  0.0169718304280
 0.7124374628501   0.0124569780990  0.0170088532421
 0.2763025250863   0.0121741311386  0.0170953520675
 0.0109658368561   0.4194306712466  0.0173888854559
 0.4289110517884   0.5599616067469  0.0174543962439
 0.4215420555115   0.0116475994785  0.0178406757287
 0.5711258590444   0.0118218313989  0.0178446863879
 0.5826868270511   0.4057889581177  0.0179046337552
 0.0130567806713   0.2725023750868  0.0181259756201
 0.0130760400964   0.7224712523233  0.0184784838882
 0.7263437062407   0.2602984019251  0.0185793564371
 0.0687230068637   0.0631417277210  0.0203217151777
 0.8652302101529   0.0720611837338  0.0213771661809
 0.0648599071037   0.8590433543910  0.0231916854098
 0.1483494943362   0.7888788352240  0.0274426710859
 0.0624359898396   0.1493935499354  0.0290301922340
 0.7871369011735   0.0656382042757  0.0294522738505
 0.0519104921610   0.5255635695605  0.0299436251629
 0.1543129927444   0.0716383926917  0.0307026948119
 0.2617842745603   0.0621479485288  0.0325263365863
 0.7667257872813   0.1658211554831  0.0327884208506
 0.2582103676627   0.6800119766139  0.0331234675192
 0.0679065925147   0.7571515437782  0.0346167526875
 0.5293578274804   0.4121503841107  0.0347081373976
 0.0666036150484   0.2612513087886  0.0347372049404
 0.0585675461899   0.3902236114535  0.0348528762454
 0.0644535360411   0.6373626559761  0.0348601561186
 0.6748138429151   0.0637583342061  0.0355471569975
 0.3914602310369   0.5503238090563  0.0360182996383
 0.6487701492307   0.2836728360263  0.0362926285843
 0.3946498220408   0.0605175522554  0.0381897702083
 0.5390137151933   0.0611990176936  0.0392252800118
 0.1627895082785   0.6861322141035  0.0482710125888
 0.6812436322641   0.1567968345899  0.0489912121566
 0.1542832878020   0.1667512624020  0.0497220833872
 0.2522727750445   0.2504803933395  0.0507065736986
 0.2547981532407   0.4994090649043  0.0509771994043
 0.1485580549194   0.5756023096087  0.0521360063667
 0.2930239606436   0.5656897354162  0.0523460874925
 0.2808991272310   0.1437921574248  0.0524440683552
 0.4820989592971   0.2518557535865  0.0527459644823
 0.5641878245444   0.1462966743153  0.0529449063728
 0.1307699644344   0.4489577586117  0.0542395594501
 0.1479692221948   0.3001174386829  0.0543470203419
 0.5638684222946   0.2813772089298  0.0547100548639
 0.4361157428790   0.4252053446420  0.0557288345913
 0.3603263935285   0.2599190004889  0.0577734264233
 0.4224188334674   0.1453238443303  0.0585393781623
 0.3719001833052   0.3780122703567  0.0609039250680
 0.2413645006928   0.3847563284940  0.0637273964449
integration degree=20  N=78:
 0.0089411337112   0.0086983293702  0.0021744545399
 0.9792622629807   0.0102644133744  0.0028987135265
 0.0105475382112   0.9785514202515  0.0030846029337
 0.0023777061947   0.0636551098604  0.0034401633104
 0.0630425115795   0.0041506347509  0.0041898472012
 0.9308422496730   0.0048053482263  0.0044738051498
 0.0629076555490   0.9316790069481  0.0047054420814
 0.9315962246381   0.0626264881801  0.0048867935750
 0.0061951689415   0.9293587058564  0.0051927643369
 0.0287125819237   0.0310202122997  0.0074073058981
 0.9293844478305   0.0342152968219  0.0079755410301
 0.0375457566621   0.9257868884669  0.0083550522910
 0.0086895739064   0.1584971251510  0.0096166660864
 0.1547597053965   0.8363606657688  0.0096318257850
 0.8331025294185   0.0089257244824  0.0098577460758
 0.8374231073526   0.1529167304078  0.0102657880301
 0.1559362505234   0.0094966240058  0.0103188103111
 0.0098599642095   0.8342211493596  0.0106291001630
 0.4055873733289   0.0074389302008  0.0106881306895
 0.5964727898618   0.3956330809311  0.0106969021010
 0.0080747800416   0.4031319425903  0.0109026461714
 0.0075073977721   0.5851609594681  0.0109899783575
 0.3936764519237   0.5974896592899  0.0113423055229
 0.5846530726212   0.0087250464968  0.0120535642930
 0.4870804112120   0.0202129229912  0.0139619193821
 0.2683512811785   0.7202340088668  0.0141147991536
 0.7223956288748   0.2662399366456  0.0141930347046
 0.2716826742357   0.0112882698808  0.0144212676268
 0.0112580842046   0.7169695963325  0.0144704346855
 0.0115034734370   0.2740067110166  0.0144949769872
 0.7140525900564   0.0113511560497  0.0145386775694
 0.4902871053112   0.4936491841468  0.0145964190926
 0.0201423425209   0.4832573459601  0.0147314578466
 0.0361107464859   0.0935679501582  0.0167463963304
 0.8607998819851   0.0397379067075  0.0168955500458
 0.1005891526001   0.8586343419352  0.0169422662884
 0.0918740717058   0.0395513001973  0.0173070172095
 0.8604888296191   0.0966224057079  0.0174524546493
 0.0439842178673   0.8561886349107  0.0177217222159
 0.2011017606735   0.7449115835626  0.0282824024023
 0.7449993726263   0.0536865638166  0.0284996712488
 0.0532186641310   0.1963754275935  0.0285005646539
 0.7453984647401   0.1982065805550  0.0300647223478
 0.1957289932876   0.0555713833156  0.0302031277082
 0.1092532057988   0.6100036182413  0.0303987136077
 0.0567625702001   0.7409121894959  0.0305668796074
 0.0483837933475   0.6075135660978  0.0306067413002
 0.1080612809760   0.1122081510437  0.0309330068201
 0.6185605900991   0.2698753703035  0.0309773820835
 0.7721296013497   0.1114117395333  0.0313146250545
 0.6115734801133   0.3389367677931  0.0313573493392
 0.3381326103376   0.0494693938787  0.0314320469287
 0.1173084128254   0.7696451309795  0.0315182143894
 0.2674551260596   0.1115718808154  0.0324248137985
 0.6542100160026   0.1906548314700  0.0347512152386
 0.0538297481158   0.3358616826849  0.0350393454927
 0.1848840324117   0.1551831523851  0.0350717420310
 0.3376267104744   0.6081402596294  0.0352129215334
 0.6067102034499   0.0542632795598  0.0352615504981
 0.4612614085496   0.0688176670722  0.0366403220343
 0.1525465365671   0.6510240845749  0.0367733107670
 0.0700582543543   0.4661904392742  0.0371675662937
 0.4704201379032   0.4634826455353  0.0373371571606
 0.1216461693746   0.2381494875516  0.0403973346588
 0.6371404052702   0.1238399384513  0.0413580040638
 0.2379904515119   0.6370216452326  0.0421957791870
 0.1483929857177   0.4894188577780  0.0495451004037
 0.3598069571550   0.1452880866253  0.0500419261141
 0.4941441055095   0.3610216383818  0.0505794587115
 0.1440630687981   0.3513508341887  0.0520037210188
 0.5019764440004   0.1435491663293  0.0521533567886
 0.3555423834298   0.5016491599502  0.0524899152358
 0.2443439540771   0.2406052129104  0.0599159762516
 0.2437064989342   0.5109017277055  0.0599609997426
 0.5122200807321   0.2452737973543  0.0599915272129
 0.2526038315178   0.3700319555094  0.0634133183449
 0.3759895652851   0.2505406611631  0.0635311861108
 0.3729077987144   0.3753750277549  0.0637206605672
integration degree=21  N=91:
 0.0035524391922   0.0035524391922  0.0006704436439
 0.0035524391922   0.9928951216156  0.0006704436439
 0.9928951216156   0.0035524391922  0.0006704436439
 0.9553548273730   0.0087898929093  0.0045472608074
 0.0358552797177   0.0087898929093  0.0045472608074
 0.9553548273730   0.0358552797177  0.0045472608074
 0.0087898929093   0.0358552797177  0.0045472608074
 0.0087898929093   0.9553548273730  0.0045472608074
 0.0358552797177   0.9553548273730  0.0045472608074
 0.8865264879047   0.1082329745017  0.0052077585320
 0.8865264879047   0.0052405375935  0.0052077585320
 0.0052405375935   0.1082329745017  0.0052077585320
 0.0052405375935   0.8865264879047  0.0052077585320
 0.1082329745017   0.8865264879047  0.0052077585320
 0.1082329745017   0.0052405375935  0.0052077585320
 0.0466397432150   0.9067205135700  0.0065435432887
 0.0466397432150   0.0466397432150  0.0065435432887
 0.9067205135700   0.0466397432150  0.0065435432887
 0.2075720456946   0.0082759241284  0.0092737841533
 0.2075720456946   0.7841520301770  0.0092737841533
 0.7841520301770   0.2075720456946  0.0092737841533
 0.0082759241284   0.7841520301770  0.0092737841533
 0.0082759241284   0.2075720456946  0.0092737841533
 0.7841520301770   0.0082759241284  0.0092737841533
 0.0858119489725   0.0314836947701  0.0095937782623
 0.8827043562574   0.0314836947701  0.0095937782623
 0.0314836947701   0.0858119489725  0.0095937782623
 0.0858119489725   0.8827043562574  0.0095937782623
 0.8827043562574   0.0858119489725  0.0095937782623
 0.0314836947701   0.8827043562574  0.0095937782623
 0.6688778233826   0.0095150760625  0.0114247809167
 0.0095150760625   0.3216071005550  0.0114247809167
 0.0095150760625   0.6688778233826  0.0114247809167
 0.6688778233826   0.3216071005550  0.0114247809167
 0.3216071005550   0.6688778233826  0.0114247809167
 0.3216071005550   0.0095150760625  0.0114247809167
 0.4379999543113   0.0099859785681  0.0117216964174
 0.0099859785681   0.5520140671206  0.0117216964174
 0.4379999543113   0.5520140671206  0.0117216964174
 0.0099859785681   0.4379999543113  0.0117216964174
 0.5520140671206   0.4379999543113  0.0117216964174
 0.5520140671206   0.0099859785681  0.0117216964174
 0.7974931072148   0.0405093994119  0.0188197155232
 0.0405093994119   0.1619974933734  0.0188197155232
 0.0405093994119   0.7974931072148  0.0188197155232
 0.1619974933734   0.7974931072148  0.0188197155232
 0.7974931072148   0.1619974933734  0.0188197155232
 0.1619974933734   0.0405093994119  0.0188197155232
 0.3864215551955   0.3864215551955  0.0235260980271
 0.3864215551955   0.2271568896090  0.0235260980271
 0.2271568896090   0.3864215551955  0.0235260980271
 0.8090129379329   0.0954935310336  0.0235571466151
 0.0954935310336   0.8090129379329  0.0235571466151
 0.0954935310336   0.0954935310336  0.0235571466151
 0.2745425238718   0.0479840480721  0.0268246207430
 0.0479840480721   0.6774734280561  0.0268246207430
 0.6774734280561   0.0479840480721  0.0268246207430
 0.6774734280561   0.2745425238718  0.0268246207430
 0.2745425238718   0.6774734280561  0.0268246207430
 0.0479840480721   0.2745425238718  0.0268246207430
 0.4053472446667   0.5429849622344  0.0314289776779
 0.0516677930989   0.4053472446667  0.0314289776779
 0.4053472446667   0.0516677930989  0.0314289776779
 0.5429849622344   0.0516677930989  0.0314289776779
 0.0516677930989   0.5429849622344  0.0314289776779
 0.5429849622344   0.4053472446667  0.0314289776779
 0.1877738615539   0.1068148267588  0.0337196192159
 0.7054113116872   0.1877738615539  0.0337196192159
 0.7054113116872   0.1068148267588  0.0337196192159
 0.1068148267588   0.7054113116872  0.0337196192159
 0.1877738615539   0.7054113116872  0.0337196192159
 0.1068148267588   0.1877738615539  0.0337196192159
 0.1195059712009   0.3057122990643  0.0427745294213
 0.1195059712009   0.5747817297348  0.0427745294213
 0.5747817297348   0.1195059712009  0.0427745294213
 0.5747817297348   0.3057122990643  0.0427745294213
 0.3057122990643   0.5747817297348  0.0427745294213
 0.3057122990643   0.1195059712009  0.0427745294213
 0.5981245743363   0.2009377128319  0.0441138932737
 0.2009377128319   0.5981245743363  0.0441138932737
 0.2009377128319   0.2009377128319  0.0441138932737
 0.2160775200005   0.3121360256673  0.0461469594684
 0.3121360256673   0.2160775200005  0.0461469594684
 0.2160775200005   0.4717864543321  0.0461469594684
 0.3121360256673   0.4717864543321  0.0461469594684
 0.4717864543321   0.3121360256673  0.0461469594684
 0.4717864543321   0.2160775200005  0.0461469594684
 0.4376579903849   0.4376579903849  0.0469152468624
 0.4376579903849   0.1246840192303  0.0469152468624
 0.1246840192303   0.4376579903849  0.0469152468624
 0.3333333333333   0.3333333333333  0.0551199980347
integration degree=23  N=105:
 0.0087809303836   0.9903676436772  0.0006438298261
 0.9903675314220   0.0087809216232  0.0006438413076
 0.0027029276450   0.0335914404439  0.0010134735710
 0.0335909214524   0.0027028946710  0.0010134752576
 0.0091675068606   0.0091676353051  0.0019679929935
 0.9675568182558   0.0084737176656  0.0033467313784
 0.0084737200688   0.9675569435345  0.0033467339208
 0.0078781948792   0.0676784943862  0.0042873323375
 0.0676785477700   0.0078781659291  0.0042873459885
 0.9470266955047   0.0442974541187  0.0043003801372
 0.0442974755680   0.9470266676487  0.0043003849098
 0.9144243214882   0.0081735455132  0.0056934629205
 0.0081735424459   0.9144244234031  0.0056934640134
 0.2497452292741   0.3833232434720  0.0061643868015
 0.3833232646055   0.2497451268005  0.0061644756418
 0.8876850353557   0.1035328809446  0.0062014513591
 0.1035329228297   0.8876849931840  0.0062014531952
 0.0077255923618   0.1403190991974  0.0069636330294
 0.1403192425107   0.0077255934624  0.0069636331842
 0.8104591009652   0.1809642523926  0.0075066257720
 0.1809643003717   0.8104590515334  0.0075066264565
 0.8330767948684   0.0083010939677  0.0079074768339
 0.0083010907126   0.8330768545392  0.0079074772485
 0.0348407706147   0.0348406969482  0.0080353344623
 0.2740287679608   0.7173981847948  0.0087963441074
 0.7173982224778   0.2740287304386  0.0087963448112
 0.2394976858234   0.0081859182262  0.0091304195716
 0.0081859185845   0.2394975566677  0.0091304213611
 0.0068836152075   0.4843740892687  0.0092821748751
 0.4843741485699   0.0068836232949  0.0092821815662
 0.4960767772741   0.4960767529507  0.0094499806178
 0.6112936776245   0.3804323691239  0.0094627468484
 0.3804323980345   0.6112936466533  0.0094627485294
 0.7303890713524   0.0083987179701  0.0095555772285
 0.0083987168639   0.7303890895407  0.0095555792843
 0.6128525675612   0.0075475979695  0.0096138842488
 0.0075475961037   0.6128525484582  0.0096138846826
 0.0079525316513   0.3559773826721  0.0099991524212
 0.3559774870460   0.0079525358502  0.0099991551850
 0.9110236977966   0.0437233665345  0.0100301319277
 0.0437233605166   0.9110236807446  0.0100301346636
 0.0388480061835   0.0967030908282  0.0124936676185
 0.0967032117936   0.0388479942386  0.0124936726125
 0.0873226911312   0.0873226620391  0.0140197309137
 0.0421445202084   0.8485617789108  0.0143336216896
 0.8485617974961   0.0421445420915  0.0143336272125
 0.8477921333864   0.1067435942472  0.0153604142740
 0.1067435889398   0.8477921328146  0.0153604183425
 0.1833966521991   0.0416340521608  0.0184523825614
 0.0416340541167   0.1833965196930  0.0184523863146
 0.7611632251560   0.1941599202852  0.0195833983573
 0.1941599254144   0.7611632153938  0.0195834019994
 0.7579378747173   0.0439826608586  0.0197632751342
 0.0439826512395   0.7579378242308  0.0197632766677
 0.0369760535918   0.5363186076436  0.0198806391019
 0.5363187134342   0.0369760780935  0.0198806485776
 0.1001256948921   0.7912267093545  0.0207181838484
 0.7912266693524   0.1001257554673  0.0207181934893
 0.0379866714177   0.4157413128558  0.0208943071440
 0.4157414028965   0.0379867061535  0.0208943251956
 0.6507106491463   0.0420141226713  0.0214864573885
 0.0420141133438   0.6507105645084  0.0214864586007
 0.0425548444254   0.2920626023484  0.0222218133036
 0.2920627107240   0.0425548546753  0.0222218160203
 0.5389729538180   0.4193031469005  0.0223345305455
 0.4193031828489   0.5389729093610  0.0223345378739
 0.6549472009700   0.3007352636162  0.0224758924946
 0.3007352790917   0.6549471812731  0.0224758980440
 0.3752400771585   0.3453980130752  0.0229701395845
 0.3453980282786   0.3752400695673  0.0229703394438
 0.0994532168761   0.1598308695187  0.0232798376102
 0.1598309359585   0.0994531960132  0.0232798427506
 0.1797326661667   0.7124585430924  0.0269483199647
 0.7124584461943   0.1797327722240  0.0269483307107
 0.1066065678636   0.7001701784175  0.0280438758010
 0.7001701904096   0.1066065855677  0.0280438764607
 0.0993303629801   0.6065647984796  0.0287526270172
 0.6065648052521   0.0993303896769  0.0287526387271
 0.1023223542704   0.2533381579528  0.0298980829063
 0.2533382324938   0.1023223826189  0.0298980922759
 0.6166226715217   0.2769502060575  0.0309004358516
 0.2769500693109   0.6166227900624  0.0309004385956
 0.0904184571873   0.4981522637001  0.0314031017088
 0.4981522767248   0.0904185045149  0.0314031073955
 0.0928231860168   0.3738418516908  0.0319191553024
 0.3738418699229   0.0928232584790  0.0319191668378
 0.2521678840407   0.2521680925697  0.0321429924062
 0.5087500218708   0.3905580544330  0.0330395601388
 0.3905579116731   0.5087501437661  0.0330395631829
 0.1706141469096   0.5266738039554  0.0356169095589
 0.5266737761312   0.1706142257537  0.0356169276054
 0.3487581527629   0.2588055084886  0.0365741189998
 0.2588053596017   0.3487583491703  0.0365741515204
 0.1696614558053   0.3013522183964  0.0365977646990
 0.3013521806875   0.1696615963219  0.0365978053889
 0.2580202409759   0.4584741774478  0.0369945680114
 0.4584740860198   0.2580203819011  0.0369945775059
 0.1848898683498   0.1848898704551  0.0374053623787
 0.6130740338465   0.1921611994069  0.0375550258317
 0.1921611750994   0.6130740398389  0.0375550312530
 0.4180541160599   0.1650613336416  0.0388887693486
 0.1650612642036   0.4180541199244  0.0388887708342
 0.5159205739625   0.2982719005229  0.0392705643548
 0.2982718935750   0.5159205534362  0.0392705802517
 0.4098894602340   0.4098894317792  0.0398766879831
integration degree=25  N=120:
 0.0082881595033   0.9848202768869  0.0014873417859
 0.4618422030241   0.5381577969759  0.0014889035262
 0.0071066441239   0.0080842361390  0.0015005944380
 0.9847613141699   0.0070015755134  0.0015059208313
 0.5374447869049   0.4625552130951  0.0015318868715
 0.0000000000000   0.4887676880140  0.0023032634487
 0.4914131929361   0.0000000000000  0.0023649067042
 0.0070345937020   0.9574158053697  0.0028751143611
 0.9564734714228   0.0364655449485  0.0029862488735
 0.0370198792045   0.0070908577166  0.0030384162737
 0.1024124542747   0.8936125594937  0.0032092459688
 0.5928065811509   0.0049451705600  0.0037029598435
 0.0050948422371   0.0996676659189  0.0037407186035
 0.0081562023689   0.0415561148784  0.0038452543223
 0.0424936107568   0.9494865260352  0.0038670778668
 0.9495543500844   0.0081794507292  0.0039192555178
 0.8932787471239   0.0053224326262  0.0039573282688
 0.0069317612927   0.9065401020433  0.0044032251724
 0.9035839030665   0.0894771171077  0.0045907108173
 0.0905665738209   0.0070525342005  0.0047023669435
 0.0083929332787   0.6663179931111  0.0050014843818
 0.6261245686071   0.0092197583153  0.0052387830156
 0.0062801592979   0.8335207460527  0.0054422104092
 0.8272539257367   0.1665134939330  0.0056931248912
 0.0062005875353   0.7424693255229  0.0059107422989
 0.1676900311185   0.0065717743528  0.0059687967687
 0.7199353069567   0.0064354534962  0.0067262190287
 0.2749740090237   0.7185296120719  0.0068307848624
 0.0079257582005   0.1766411374714  0.0069531259112
 0.0069981220752   0.2704767254004  0.0072460270642
 0.8125248773263   0.0082299533210  0.0072728189613
 0.0073536969970   0.5934167875453  0.0073008930847
 0.7283665935411   0.2648817553752  0.0073604666776
 0.1800642304565   0.8115848976682  0.0074119923255
 0.2658102467762   0.0068553525429  0.0074892214336
 0.0070892364520   0.3757632659744  0.0078604067260
 0.3774054302043   0.6148573533757  0.0078621726423
 0.0369649608668   0.9210792302893  0.0080506361066
 0.9203194109805   0.0426025082114  0.0081442860473
 0.0425477806431   0.0372689941794  0.0081478804152
 0.6191278394983   0.3724055713809  0.0092444146612
 0.3762697209178   0.0081436422011  0.0094674635165
 0.0956111149690   0.8771098372601  0.0097132210137
 0.0302473410377   0.0943858903393  0.0099753581151
 0.8739905691754   0.0313198990883  0.0103367803673
 0.8604133734958   0.1049019782046  0.0112263277166
 0.0347307852352   0.8609856462886  0.0114309118745
 0.1043606608343   0.0357152881004  0.0115550567487
 0.7797622824754   0.1872318199265  0.0135575856957
 0.0185865164256   0.4834397678794  0.0135984962900
 0.0324585286618   0.7783474916042  0.0137754813837
 0.8371293901157   0.0804060570156  0.0137961015942
 0.0836602075315   0.8421414817051  0.0138408839904
 0.0784070242501   0.0849927089145  0.0140634019977
 0.4929238648458   0.4892855914710  0.0140991451009
 0.1870637584073   0.0345210858281  0.0142004111991
 0.4892636967025   0.0190774755077  0.0144518424517
 0.0401982618372   0.1691143187109  0.0150245979639
 0.7894259278865   0.0412206731484  0.0152817804122
 0.1686260456429   0.7894860640585  0.0155550724169
 0.3750901913174   0.5895318272013  0.0164570886000
 0.0356362876880   0.3681256217699  0.0165275759573
 0.5887548164804   0.0359968962541  0.0166847554451
 0.0373308082182   0.6790704673533  0.0167409312985
 0.2820769993374   0.0373639992361  0.0168674663361
 0.6819277603320   0.2803330345725  0.0168882230165
 0.0374938324382   0.2634016180014  0.0172087112691
 0.6984079204127   0.0364154673322  0.0174681068264
 0.2654390894079   0.6980717436193  0.0176663899614
 0.1429848440800   0.7612254618453  0.0182967621475
 0.7623554007647   0.0943741220275  0.0183576852459
 0.0934222022749   0.1479799836832  0.0186392569521
 0.5759004479923   0.3821329641698  0.0189781060590
 0.3822427332525   0.0426716362301  0.0191847922578
 0.0411414081675   0.5718082874432  0.0194080442044
 0.0802462538379   0.7702204382042  0.0194720072193
 0.7625229819410   0.1559420577362  0.0200855080495
 0.1524941445131   0.0842965421322  0.0201673909332
 0.0622159195833   0.4538181318873  0.0221742162761
 0.1109539036076   0.4586014071171  0.0229702440508
 0.4575627212057   0.4795313560210  0.0233465117399
 0.4322865136374   0.1230591237472  0.0234883135338
 0.5865002850241   0.0834119779793  0.0240682099018
 0.0869359250818   0.6755677013351  0.0240910792953
 0.0929594906936   0.2326500892727  0.0245677049481
 0.6661932141454   0.2448294007406  0.0246536315719
 0.4780306362227   0.0661749044835  0.0246756530052
 0.4372215294577   0.4442145585244  0.0249704602710
 0.6779224504669   0.0929096534577  0.0250026544082
 0.2423431255660   0.0889793655129  0.0250490869426
 0.2288925420305   0.6780053081672  0.0250936250125
 0.3315065049959   0.5847381559741  0.0251482076226
 0.3424200526607   0.5139245722736  0.0255010290447
 0.0862630046475   0.3340976249234  0.0256544511979
 0.5113188946635   0.1380154720554  0.0257974750630
 0.1538977841001   0.6788062619562  0.0270007753993
 0.6779951348472   0.1663358925269  0.0274431536844
 0.1664600469411   0.1582214504849  0.0277072401488
 0.0950910318888   0.5666590332543  0.0278284415364
 0.3436048136712   0.0978960873457  0.0287207381105
 0.5560417025366   0.3468917820947  0.0288826834956
 0.1452404029513   0.3599534491052  0.0293302729759
 0.1619685156238   0.5810131373330  0.0318902879557
 0.5800164844262   0.2560674640672  0.0319083660286
 0.2450201223288   0.5881469552102  0.0320938960329
 0.2557621891794   0.1652244065047  0.0321618608780
 0.2205239985511   0.3496507466106  0.0322424127534
 0.4940183111285   0.2549448448453  0.0327072446421
 0.2531570689798   0.2543369115017  0.0329946316695
 0.5846891116357   0.1666603916479  0.0331828096025
 0.1660333602278   0.2523240191705  0.0334857162651
 0.2505426292461   0.4959007627528  0.0335468472792
 0.3519336802182   0.1805380367800  0.0337049042988
 0.3502668835419   0.4358582329881  0.0340361462767
 0.4400892485512   0.2120576104941  0.0342465235323
 0.4680855471546   0.3552681570774  0.0345528817251
 0.1770237763947   0.4670352922266  0.0356782875703
 0.3900920779501   0.3323152819300  0.0364656225016
 0.2805847774120   0.3898041176680  0.0365172708706
 0.3361523347440   0.2778500044356  0.0371924811018